\documentclass[12pt,reqno,a4paper]{amsart} 
\usepackage{amsthm}
\usepackage{amsmath}
\usepackage{amsfonts}
\usepackage{mathtools} 			
\usepackage[pdfpagelabels,hidelinks]{hyperref}
\usepackage{enumitem}
\usepackage{bbm}
\usepackage{stmaryrd}
\usepackage{lmodern}
\usepackage{kpfonts}
\numberwithin{equation}{section}
\usepackage[top=1.5cm, bottom=1.5cm, left=1.5cm,right=1.5cm]{geometry}
\usepackage{stmaryrd}
\usepackage{mathrsfs}
\numberwithin{equation}{section}

\usepackage{graphicx} 
\usepackage[font=scriptsize,justification=centering]{caption}
\usepackage{subcaption}
\newtheorem{thm}{Theorem}[section]
\newtheorem{cor}[thm]{Corollary}

\newtheorem{ppn}[thm]{Proposition}

\theoremstyle{definition}

\newtheorem{remark}[thm]{Remark}

\usepackage{xcolor}

\definecolor{nahu}{rgb}{0.8, 0.5, 0.2}

\newcommand{\CAL}[1]{\mathcal{#1}}  
\newcommand{\BB}[1]{\mathbb{#1}}
\newcommand{\tn}[1]{\textnormal{#1}}

\DeclareMathOperator{\dd}{\textnormal{d}\!}
\def\nn{\nonumber}
\def\la{\langle}
\def\ra{\rangle}
\def\R{\mathbb R}
\def\N{\mathbb N}
\def\X{\mathcal X}
\def\L{\mathcal L}
\def\I{\mathbb I}

\def\maB{\mathcal B}
\newcommand\pr[1]{{\mathbb P}\left(#1\right)}

\DeclareMathOperator{\bb1}{\mathbbm{1}}
\DeclareMathOperator{\ncond}{\textnormal{\textsc{Ncond}}}
\DeclareMathOperator*{\argmax}{arg\,max}

\DeclareMathOperator{\supp}{\textnormal{supp}}

\allowdisplaybreaks
\title{Generalized max-weight policies in stochastic matching}
\author[M. Jonckheere]{M. Jonckheere}
\address{Instituto de C\'alculo, Universidad de Buenos Aires, Buenos Aires, Argentina}
\email{matthieu.jonckheere@gmail.com}
\author[P. Moyal]{P. Moyal}
\address{Institut Elie Cartan, Université de Lorraine}
\email{pascal.moyal@univ-lorraine.fr}
\author[C. Ram\'irez]{C. Ram\'irez}
\address{Aristas}
\email{c.ramirez@aristas.com.ar}
\author[N. Soprano-Loto]{N. Soprano-Loto}
\address{Facultad de Ingenier\'ia, Universidad de Buenos Aires, Buenos Aires, Argentina and Aristas}
\email{nsoprano@fi.uba.ar}

\begin{document}

%\large
\begin{abstract}
We consider a matching system where items arrive one by one at each node of a compatibility network according to Poisson processes and depart from it as soon as they are matched to a compatible item. The matching policy considered is a generalized max-weight policy where decisions can be noisy. Additionally, some of the nodes may have impatience, i.e. leave the system before being matched.
Using specific properties of the max-weight policy, we construct several Lyapunov functions, including a simple quadratic one. 
This allows us to establish stability results,
to construct bounds for the stationary mean and variances of the total amount of customers in the system, and to prove exponential convergence speed towards the stationary measure.
We finally illustrate some of these results using simulations on toy examples.
\end{abstract}

\maketitle

%\tableofcontents

%\large

\section{Introduction}
The so-called {\em bipartite stochastic matching model} was introduced in \cite{CKW09} (see also \cite{AdWe}) as a variant of skill-based service systems.
In this model,
customer/server couples enter the system at each time point.
A matching occurs when an arriving customer (resp. server)
finds a compatible server (resp. customer) present in the system,
in which case both leave instantaneously.
Otherwise, items remain in the system waiting for  compatible arrivals
(in particular, there are no service times).
The field of applications of such models is very wide, from housing to job applications, from transplant protocols to blood banks, peer-to-peer systems to dating websites, ride sharing, etc.
The mathematical setting is the following:
there is a set $C$ of customer classes and a set $S$ of server classes, and a bipartite graph on the bipartition $C\cup S$ indicates the compatibilities:
a customer class $c\in C$ and a server class $s \in S$ are compatible if and only if there is an edge between $c$ and $s$ in the compatibility graph.
Customer/server couples  enter the system in discrete time, and the classes of the incoming couples are i.i.d. from the measure $\mu:=\mu_C\otimes\mu_S$ on $C\times S$.
(In particular, importantly, the class of the incoming customer and server are independent of one another.)
%Each incoming customer (resp.server) is matched with a compatible server (resp. customer) present in the system unless there is none, in which case it is stored in a buffer.
In case of multiple compatibilities,
it is the role of the {\em matching policy} to determine what match is performed.
In the two aforementioned references, the matching policy is `First Come, First Served'. 
In \cite{ABMW17} it is shown that, under a natural stability condition, the stationary probability of the latter system has a remarkable product form. Interestingly, this result can then be adapted to 
various skill-based queueing models as well, and in particular those applying (various declinations of) the so-called FCFM-ALIS (Assign the Longest Idling Server) service discipline ---see e.g. \cite{AW14}. 
For various extensions of bipartite matching models, see also \cite{AKRW18}. \cite{BGM13} then extends the above settings to the case where the classes of the incoming couples are not drawn independently, and to a more general class of matching policies.
These models are then called `Extended Bipartite Models'.
In particular, the stability condition in \cite{ABMW17} is shown to be also sufficient for the stability of extended bipartite models ruled by the `Match the Longest' policy, consisting in always matching incoming items to the compatible class option
having the longest queue.
In \cite{MBM18}, a Coupling From The Past result is obtained, showing the existence of a unique bi-infinite matching in various cases of extended bipartite matching models, and for a broader class of matching policies than FCFS, thereby generalizing various results of \cite{ABMW17}. 

In many practical systems, it is often more realistic to assume that arrivals are simple instead of pairwise. 
Moreover, the bipartition of the compatibility graph into classes of {\em servers} and {\em customers} only suits applications in which this bipartition is natural (donor/receiver, house/applicant, job/applicant, and so on).
However, in many cases, the context requires that the compatibility graph take a general (i.e., not necessarily bipartite) form. For instance, in dating websites, it is a priori not possible to split items into two sets of classes with no possible matches within those sets. Similarly, in kidney exchange programs, intra-incompatible  donor/receiver couples enter the system looking for a compatible couple to perform a `crossed' transplant. 
Then, it is convenient to represent donor/receiver couples as  {\em single} items, and compatibility between couples 
means that a kidney exchange can be performed between the two couples (the donor of the first couple can give to the receiver of the second, and the donor of the second can give to the receiver of the first). In particular, if one considers blood types as a primary compatibility criterion, the compatibility graph between couples is naturally non-bipartite. Motivated by such applications, among others, 
 \cite{MaiMoy16} proposed a generalization of the above matching models, termed {\em general stochastic matching model}, in which items enter one by one, and the compatibility graph $G$ is general. 
For a given compatibility graph $G$ and a given matching policy $\Phi$, the {\em stability region} {\sc Stab}$(G,\Phi)$ of the latter system is defined as the set of those measures $\mu$ deciding the class of the incoming items, such that the Markov chain representing the system is positive recurrent. \cite{MaiMoy16} identified an universal necessary condition for stability, that is, a set {\sc Ncond}$(G)$ 
that includes {\sc Stab}$(G,\Phi)$ for any matching policy $\phi$ and any graph $G$. Following \cite{BGM13}, \cite{MaiMoy16} also showed that 
`Match the Longest' is {\em maximal stable}, i.e. it has a maximal stability region, as the two sets {\sc Ncond}$(G)$ and {\sc Stab}$(G,\Phi)$ coincide. \cite{MoyPer17} then showed that, aside for a particular class of graphs, the random `Class-uniform' policy (choosing the match in a compatible class of the incoming item, uniformly at random among those having a non-empty queue) is never maximal-stable, and that priority policies are in general {\em not} maximal-stable. Then, \cite{MBM17} showed that the policy `First Come, First Matched' is maximal-stable, and that, similarly to \cite{ABMW17}, the stationary probability of the system can be written in a product form. Various extensions of the general stochastic matching model can also be found in \cite{BC15,BC17} regarding systems with bipartite complete compatibility graphs, but in which each match can be enacted with a nominal probability, and in \cite{RM19} to stochastic matching systems on hypergraphs, thereby matching items by groups of two or more. 
In a recent line of study, stochastic matching models have also been studied from the point of view of stochastic optimization, see \cite{BM14,GW14,NS19}. 
In particular, the Max-Weight policy introduced in \cite{MW} is a natural policy allowing a natural trade-off between efficiency and stability.

Here we consider a matching system where items arrive one by one at each node of a compatibility network according to Poisson processes and depart from it as soon as they are matched to a compatible item. The matching policy considered is a generalized max-weight policy where decisions can be subject to noise. Additionally, some of the nodes may have impatience, i.e. leave the system before being matched.

\subsection*{Contributions}
In many applications, reneging of items need to be taken into account: for instance, applicants may renege from the systems, whereas jobs/houses/university spots may no longer be available. 
Similarly, in organ transplants, patients may renege from the system due to death or change in their diagnosis, whereas organs have a short `life duration' outside a body, and need to be transplanted within a very short period of time ---see a first approach of such systems in \cite{BDPS11}.

The present work is dedicated to an extension of general stochastic matching models to the case where (i) several classes of items are impatient, in that they 
renege from the system if they do not find a match within a given (random) period of time and (ii) the matching policy is of the Max-Weight class, with possible errors in the 
assessment of the system state. Specifically, we consider that the policy decisions may be subject to noise that can lead to errors in the estimate of the queue lengths, which is necessary to implement 
matching policies of the Max-Weight type. 
This situation has been considered for systems of parallel queues in \cite{MP19}.
%and it is the spirit of the `Power of $ d $' policies (see \cite{Mit96} and related references), to reduce the complexity of routing and the probability of errors, by routing the requests only to a subset of the available servers.}

Using the specific structure of the Max-Weight policy, our first contribution is to prove that a simple quadratic function is 
a Lyapunov function when a traffic condition, that can be understood as a generalization to the case of reneging of the condition ``$\mu\in\textsc{Ncond}(G)$'', is satisfied (see (\ref{eq:Ncond}) below).  
This allows us directly to characterize the stability region and, as a corollary,
to conclude that the Max-Weight policy is maximum stable.
It is worth mentioning that the last one is a new result even in the case without impatience and noise.
%retrieve the known stability results in the case without impatience and noise, and to generalize the stability results in this more general context.
The proof of the Lyapunov property is quite intricate
and relies on carefully exploiting this $\ncond$ condition.
We can then construct several additional Lyapunov functions satisfying interesting drift inequalities.
When combined with stochastic comparisons, this allows us to  construct bounds for the stationary mean and variances of the total amount of customers in the system, and additionally  to prove exponential convergence speed towards the stationary measure.

We finally illustrate some of those results using simulations on toy examples.

\subsection*{Outline of the paper.}
In Section \ref{section:the_model},
we give the formal definition of the model we study.
Section \ref{section:results} is dedicated to the theoretical results.
We start in Subsection \ref{subsection:quadratic} with 
Proposition \ref{prop:quadratic},
a key result in which we control the drift of the quadratic function.
Subsection \ref{subsection:stability} is dedicated to stability:
a stability characterization is given in Theorem \ref{thm:stability},
and Proposition \ref{prop:Ncond} describes precisely whether the model is stabilizable.
A geometric Lyapunov bound is given in Proposition \ref{ppn:geometric_bound},
and the exponentially fast convergence to equilibrium is stated as a consequence in Corollary 
\ref{cor:exponential_convergence}.
In Subsection \ref{subsection:bounds}, some bounds for the size of the largest queue in the stationary regime are presented:
a lower bound for the expectation in Proposition \ref{ppn:lower}, and upper bounds for the expectation and the variance in
Corollaries \ref{cor-1st_moment} and \ref{cor:variance}.
Section \ref{section:simulations} is dedicated to simulations:
we compare the Max-Weight policy with a priority policy in Subsection \ref{subsection:greedy}, 
and an analysis of the theoretical bounds obtained in the previous section is given in Subsection \ref{subsection:bounds}.
Finally, all the proofs are in Section \ref{section:proofs}.

\section{The model}\label{section:the_model}

\subsection{Notation} 
The graph $ G=(V,E) $ is assumed to be finite, simple, undirected and connected.
The non-connected case can be achieved by decomposing in connected components.
Any edge $i-j$ in $E$, for $i,j\in V$, is denoted by $\la i,j \ra$.
We say that a subset $ I \subset V$ is an independent set if $ I\neq \emptyset $ and if $ i,j\in I $ implies $ \la i,j\ra \notin E $.
%In words, all the clusters  of  $I$ are singletons.
For any subset $A \subset V$, we let 
\begin{equation}
\label{eq:defIA}
\BB I(A)=\{ A'\subset A: A'\mbox{ is an independent set} \}
\end{equation}
\noindent in a way that $\I:=\I(V)$ gathers all independent sets of the graph $G$. 
For a subset $ A\subset V $, let $ E(A)=\{ i\in V: \la i,j \ra\in E\mbox{ for some }j\in A \} $.
For singletons, we write $ E(i) $ instead of $ E(\{i\}) $.
Observe that $ A\cap E(A) $ is the union of the clusters of $ A $ that are not singletons;  in particular, $ A\cap E(A)=\emptyset $ if $ A \in \BB I$.
Let $\R, \R_+,\R_+^*,\mathbb Z,\N_0$ and $\N$ denote respectively the sets of real numbers, non-negative real numbers, positive real numbers, integers, non-negative integers and positive integers. 
For $p,q\in\mathbb Z$, we let $\llbracket p,q \rrbracket:=[p,q]\cap \mathbb Z$ be the set of all integers in the closed interval $[p,q]$. 
Let also for all $q\in\N$ and for all $x=(x(1)\cdots x(q))\in\R^q$,
$\| x \|=\left(\sum_{i\in V}x(i)^2\right)^{1/2}$,
$ \|x\|_\infty =\max_{i\in V}|x(i)|$ and $ \|x\|_1=\sum_{i\in V}|x(i)|$.  

\subsection{A stochastic matching model with reneging} 
The (continuous-time) 
stochastic matching model with reneging $(G,\textsc{mw},\lambda,\gamma)$ is formally defined as follows:  we give ourselves a set of $|V|$ independent Poisson processes $N(i),\,i\in V$, of respective intensities 
$\lambda(i)>0\,,\,i\in V$, 
and  denote by $\lambda:=(\lambda(i)\,,\,i\in V)$ the {\em arrival intensity vector} of the model.
For any $i\in V$ we denote by $\{T^i_n\,:\,n\in \N\}$ the set of points of $N(i)$. For any subset $A\subset V$, denote by $\lambda(A)=\sum_{i\in A}\lambda(i)$, the intensity of the 
superposed Poisson processes having indexes in $A$. The overall superposed Poisson process of the $N(i)$'s, $i\in V$, is denoted by $N$. 
It is thus a Poisson process of 
intensity $\lambda(V)$, and its points are denoted by $\{T_n\,:\,n\in \N\}:=\bigcup_{i\in V}\{T^i_n\,:\,n\in \N\}$. 
We consider that for all $i$, upon each point $T^i_n$, $n\in\N$, an {\em item} of {\em class} $i$ (or {\em $i$-item}, for short) enters the system. 
Thus the points of $N$ are the overall arrival times in the system. We say that two items of respective classes $i$ and $j$ in $V$ are 
{\em compatible} if and only if $\la i,j \ra \in E$. To each class $i\in V$ is also associated a parameter $\gamma(i) \ge 0$.
Denote $\gamma:=(\gamma(i),\,i\in V)$, the {\em reneging rate vector} of the model. 
We then associate to each $i$-item an Exponential r.v. of distribution $\mathcal E(\gamma(i))$, interpreted as its {\em patience time}. 
A r.v. of distribution $\mathcal E(0)$ is interpreted as a.s. infinite, and we denote by 
\[R:=\{i \in V\,:\,\gamma(i)>0\}\quad\mbox{ and }\quad R^c:=\{i\in V\,:\,\gamma(i)=0\}\]
the subsets of nodes whose items do (resp., do {not}) renege from the system. 
The corresponding item then reneges from the system at the end of this time period, provided that it has not left before, 
following a matching mechanism that is defined below. 
We assume that the patience times of all items are independent of one another, and of the other r.v.'s involved. 

The dynamics of the model $(G,\textsc{mw},\lambda,\gamma)$ is then defined recursively at the successive points of $N$, as follows: 
\begin{itemize}
\item We start at time 0 with a system containing in its buffer, incompatible items. In other words, for any $\la i,j \ra \in E$ in the system at time 0 
either there are no $i$-items or no $j$-items, or neither. Such a buffer content is said {\em admissible}. 
To each of these initial items is associated a patience time that is drawn, independently of everything else, from the distribution corresponding to its class, e.g., from $\mathcal E$($\gamma(j)$) for a $j$-item. 
As previously mentioned, if $j\in R^c$, i.e., $\gamma(j)=0$, the patience time of the corresponding item is set as infinite. 
\item Let $n\in\N$, and suppose that the buffer content of the system was admissible at time $T_{n-1}$. 
We first update the buffer content at each reneging time in $[T_{n-1},T_n)$, by discarding all items present in the system at time $T_{n-1}$ and whose patience time has elapses during that time interval. Observe that the buffer content remains admissible after any of these reneging. Then, suppose that $T_n$ is a point of the $j$-th arrival process, i.e. $T_n=T^j_k$ for some $k\in \llbracket 1,n \rrbracket$, or in other words, a $j$-item enters the system at $T_n$. This happens with probability 
\begin{equation}
\label{eq:defmu}
\mu_\lambda(j):=\frac{\lambda(j)}{\lambda(V)}\cdot 
\end{equation}
We then face the following alternative:  
\begin{itemize}
\item Either the buffer contains no remaining compatible item with the incoming $j$-item, i.e. there are no stored items of classes 
in $E(j)$, in which case the latter item is stored in the buffer. We then draw the (possibly infinite) patience time of this item from the distribution $\mathcal E$($\gamma(j)$), independently of everything else.
% \nahu{[I don't understand why you treat the reneging of this item separately. Is it not included in the previous phrase ``We first update the buffer content at each reneging time in $[T_{n-1},T_n)$, by discarding all items present in the system at time $T_{n-1}$ and whose patience time has elapses during that time interval.''?]}\pasc{[The only thing we say here is that we draw the patience time of an item upon its arrival, only if it is not matched right away.]}
\item Or, there is at least one remaining item having class in $E(j)$ in the buffer. Then the incoming $j$-item must be matched 
with one of these compatible items, and it is the role of the matching policy to determine its match. 
We assume that the matching policy is of the {\em generalized Max-Weight policy} type with parameters $w$ and $F$ ({\sc mw}$_{w,F}$, or simply {\sc mw}, for short). 
To define this, let $w:=\{w_{(j,i)}\,:\,(i,j)\in V^2\mbox{ s.t. }\la i,j \ra \in E\}$ be a family of non-negative numbers and $\{F_{(j,i)}:\,(i,j)\in V^2 \}$ be a family of cumulative distribution functions (cdf's, for short). Denote, for all $i\in V$, by $x(i)$ the number of stored $i$-items at this time. Then the match of the incoming $j$-item is chosen uniformly at random from the set
\begin{align}\nn
\argmax_{i\in E(j)\,:\, x(i)>0} \big([x(i)+ U_{n,(j,i)}]^++w_{( j,i)}\big),
\end{align}
where $U_{n,(j,i)}$ is a r.v. drawn at time $n$ independently of everything else from the cdf $F_{(j,i)}$, interpreted as the measurement error of the queue length of node $i$, made upon an arrival of a $j$ element at time $T_n$, and $w_{(j,i)}$ 
is interpreted as the specific revenue of matching an incoming $j$-item with a $i$-item. (Notice that symmetry between $i$ and $j$ is a priori not assumed in both definitions.) In words, the incoming $j$-item is matched with an element of a class $i\in E(j)$ that maximizes the reward of the match 
of $j$ with $i$, where the reward is measured as a linear function of the queue sizes of the neighboring nodes of $j$, measured with an error given by $U_{n(j,i)}$ (we prevent absurd measurement of negative queue lengths by taking the positive part of the latter), plus the specific revenue 
	for matching the $j$-item with items of its neighboring classes. 
	Then, once the matching is done, the incoming $j$-item leaves the system right away, together with its match. 
%(that is, at time $T_{n}$ and just after reneging all the items whose patience has elapsed between $T_{n-1} and $T_n$) 
\end{itemize}
It is then easily seen that, in any case, the buffer content at time $T_n$ is also admissible. 
\item We continue the construction inductively. 
\end{itemize}
\begin{remark}
Particular cases of generalized Max-Weight policies are: 
\begin{enumerate}
\item Classical {\bf Max-Weight} policies as defined in cite{?}, taking the cdf $F_{(j,i)}:=\bb 1_{\R_+}$ as the one corresponding to the Dirac $\delta_0$ distribution, i.e., no error, for all $(i,j)\in V^2$. 
\item The {\bf Match the Longest} policy introduced in \cite{MaiMoy16} taking $F_{(j,i)}:=\bb 1_{\R_+}$ for all $(j,i)$ (i.e., no errors), 
$w_{(j,i)}:=w_j$ for all $j\in V$ and all $i\in E(j)$ (i.e., same revenue for all matchings with a $j$-item). 
%\item The {\bf Match the Shortest} policy also introduced in \cite{MaiMoy16} taking parameters as above, but for $\beta=-1$.
%\item Any {\bf Strict priority} policy (again, see \cite{MaiMoy16}), taking $\beta=0$. 
\end{enumerate}
\end{remark}

\subsection{Markov representation}
Let for all $t\ge 0$, $X_t:=\left(X_t(i),\,i\in V\right)$ be the queue length vector at time $t$, namely, for any $i\in V$, $X_t(i)$ is the number 
of $i$-items in the buffer at time $t$. It is easily seen that the process $X= \{X_t\,:\,t\ge 0\} $ is a Continuous-Time Markov Chain (CTMC) on 
the commutative state-space
\begin{align}\nn
\CAL X=\big\{ x\in\BB N_0^{|V|}: x(i)x(j)=0 \ \mbox{ for any } (i,j)\in V^2 \mbox{ such that }\langle i,j\rangle\in E
\big\}.
\end{align}

For any $x\in \X$, let $ \supp(x)=\{ i\in V:x(i)>0 \} $ be its support. 
Observe that $ \supp(x)\in \BB I $
for every $ x\in\X\setminus\{\mathbf 0\}$.  % (according to our definition, $ \emptyset $ is not an independent set).
%\nahu{[Apparently it is more natural to include the emptyset in the definition of an independent set, also because of a previous comment I made about $ \BB I(A) $.]}
Denote for all $i\in V$ by $\delta_i$, the $i$-th canonical vector of $\N_0^{|V|}$. 
The infinitesimal generator $\L$ of $X$ is then simply given by 
\begin{align}
\label{eq:genX}
&\L f(x)=&\sum_{i  \notin  E(\supp(x))}
\left[f(x+\delta_i)-f(x)\right]\lambda(i)
+\sum_{i  \in  \supp(x)} \left[f(x-\delta_i)-f(x)\right]\bigg[\gamma(i)x(i)+\sum_{j \in   E(i)}\lambda(j)\nu_{x,j}(i) \bigg]
\end{align}
for all $x\in \X$ and $ f:\X\to\R $, where for all $i\in V$ and $j\in E(i)$, $ \nu_{x,j}(i) $ is the probability of performing the match $ \la i,j\ra $, given that a $j$-item enters the system in state $x$.
From the very definition of {\sc mw}, the latter is fully determined by the cdf's $F_{(j,i)}$, $i\in \supp(x)\cap E(j)$.

\section{Results}\label{section:results}

\subsection{Quadratic Lyapunov function}
\label{subsection:quadratic}

Fix a matching model $(G,\textsc{mw},\lambda,\gamma)$. Recalling (\ref{eq:defIA}), we 
let % of the probability measure $\mu$ on $V$, define 
\begin{equation}\label{eq:defeta}
\eta:=\eta(\lambda)=\begin{cases}
\min\{ \lambda\left(E(I)\right)-\lambda(I) :I\in\I(R^c)\}&\mbox{ if }R^c \ne \emptyset;\\
0&\mbox{ else.} \end{cases}
\end{equation}

Let $\{U_{(j,i)}:\,(i,j) \in V^2 \mbox{ such that }\la i,j \ra \in E\}$
be a mutually independent family of r.v.'s such that, for any such $(i,j)$, 
$U_{(j,i)}$ has cdf $F_{(j,i)}$.
For any such $(i,j)\in V^2$ and any $ \kappa\in (0,\lambda(V)) $, 
define
\begin{align}\nn
u_{\kappa,(j,i)} =\inf\bigg\{u\in \R_+\,:\,\BB P\big(|U_{(j,i)}|>u\big)
<1-\Big(1-{\frac{\kappa}{\lambda(V)}  }\Big)^{\frac{1}{2|E|}}\bigg\},
\end{align}
and let
\begin{align}
u_\kappa = \max\bigl\{u_{\kappa,(j,i)}\,:\,(i,j)\in V^2\mbox{ s.t. }\la j,i \ra  \in E\bigl\}.\label{eq:defu}
\end{align}
Observe that, from the iid assumptions, for all $n$ we have 
\begin{align}\label{eq:tail}
\pr{\mathcal B_{\kappa,n}}
:=
\BB P\Big[\bigcap
\big\{|U_{n,(j,i)}| \le u_\kappa\big\}\Big]
=
\BB P\Big[\bigcap
\big\{|U_{(j,i)}| \le u_\kappa\big\}\Big]
=:\pr{\mathcal B_\kappa}\ge 1 - { \frac{\kappa}{\lambda(V)}  },
\end{align}
where the intersections are taken over the pairs $ (i,j) \in V^2$ such that $\la i,j \ra  \in E $.
The key result of our analysis is the following generator bound, 
\begin{ppn}\label{prop:quadratic}
	Define the mapping $ f_2:\CAL X\to\R $ by $ f_2(x)=\sum_{i\in V}x(i)^2 $ for all $x\in\X$. 
	Then, for all $x\in\X $ and all $\kappa \in (0,\lambda(V))$, we have that 
	\begin{align}\label{inequality-225}
	\L f_2(x)\le \lambda(V)+ 2 \sum_{i\in R}
	\left[-\gamma(i)x(i)^2+\left(\frac{1}{2} \gamma(i) +\lambda(i)\right)x(i)\right]
	+2\left[ \lambda(V)(2u_\kappa+\check w)|V|
	+(\kappa-\eta)\| x \|_{\infty}\right]\bb 1_{R^c \ne \emptyset},
	\end{align}
	for $\check w=\max_{e\in E} w_e$, $\eta$ defined by (\ref{eq:defeta}) and $u_\kappa$ defined by (\ref{eq:defu}). 
\end{ppn}

\subsection{Stability results}
\label{subsection:stability}

We say that the matching model $(G,\textsc{mw},\lambda,\gamma)$ is stable if the CTMC $ (X_{t})_{t\ge 0} $ is positive recurrent.
The next result provides a necessary and sufficient condition for stability.

\begin{thm}\label{thm:stability}
For any connected graph $G$ and any $\gamma \in (\R_+)^{|V|}$, the model $(G,\textsc{mw}_{w,F},\lambda,\gamma)$ is stable if and only if its intensity vector $ \lambda $ 
belongs to the set
\begin{align}\label{eq:Ncond} 
\ncond(G,\gamma)=\left\{\lambda\in(\BB R_+^*)^{|V|}\,:\,\lambda(I)<\lambda\left(E(I)\right) \mbox{ for any } I\in\I(R^c)\right\}.
\end{align} 
\end{thm}

\begin{remark}[Full reneging case]
In the case where $R^c=\emptyset$, i.e. all item classes  renege, then condition (\ref{eq:Ncond}) is empty, and we retrieve the good-sense result that 
any such matching model is stable, whatever the intensity vector $\lambda$. 
%Moreover, from (\ref{inequality-225}) we get that
% \[ \L f_2(x)\le 2 \sum_{i\in V}\left(-\gamma(i)x(i)^2+\left(\frac{1}{2} \gamma(i) +\lambda(i)\right)x(i)\right).\]
\end{remark}

\begin{remark}[No-reneging case]
In the case where $R=\emptyset$ (no reneging), the model amounts to a continuous-time version of that introduced in 
\cite{MaiMoy16}. We can then complete the results in \cite{MaiMoy16}, \cite{MoyPer17} and \cite{MBM17}, by observing that any model $(G,\textsc{mw},\lambda,\gamma)$ is stable whenever 
$\lambda$ belongs to the set 
\begin{align}
\label{eq:Ncondno}
\textsc{Ncond}(G)=\left\{\lambda\in(\BB R_+^*)^{|V|}\,:\,\lambda(I)<\lambda\left(E(I)\right), \mbox{ for all } I\in\I\right\}.
\end{align}
Any policy of the {\sc mw} type is then maximal, in the sense that it has the largest stability region. We can specialize this result to the case where the matching policy is 
`Match the Longest', i.e., $F_{(j,i)}:=\bb 1_{\R_+}$ for all $(j,i)$  (no errors) and $w_{(j,i)}:=w_j$ for all $j\in V$ and all $i\in E(j)$ (equal rewards for all matches involving an incoming $j$-item, $j\in V$), thereby retrieving assertion (15) in Theorem 2 of \cite{MaiMoy16}. 
\end{remark}

As can be easily seen in the proof of the necessity statement in Theorem \ref{thm:stability}, for a connected graph $G$ and a reneging rate vector 
$\gamma$, $\lambda$ being an element of $\textsc{Ncond}(G,\gamma)$ is necessary for a stable matching model defined on $G$ to exist, whatever the matching policy is. 
The sufficiency statement means that such models always exist whenever $\lambda$ do belong to $\textsc{Ncond}(G,\gamma)$ ---it suffices to implement a matching policy in the Max-Weight class. Therefore it is natural to say that the couple
$(G,\gamma)$ is {\em stabilizable} if there exists an intensity vector $\lambda$ such that the model with reneging $(G,\textsc{mw},\lambda,\gamma)$ is stable, that is, if \textsc{Ncond}$(G,\gamma)$ is non-empty. The following proposition makes precise the class of stabilizable graphs.

\begin{ppn}
\label{prop:Ncond}
Let $G=(V,E)$ be a connected graph and $\gamma\in(\R_+)^{|V|}$ be a reneging rate vector. Then \break $\textsc{Ncond}(G,\gamma)$ is non-empty if and only if 
$G$ is non-bipartite or $R\ne \emptyset$. 
\end{ppn}

\subsection{Geometric Lyapunov function}

Building on the previous results, and following classical approaches for Markov processes with bounded jumps, we can show how to define a geometric Lyapunov function.
If $ \ncond(G,\gamma) $ is satisfied,
we call $ \pi $ the unique invariant distribution whose existence is guaranteed by Theorem \ref{thm:stability}.

\begin{ppn}\label{ppn:geometric_bound}
Under the $\ncond$ assumption,
if $\alpha>0$ is small enough,
there exist $C,K >0$ such that
$$ \L e^{\alpha\|x\|} \le - C e^{\alpha\|x\|}$$
for every $x\in\CAL X$ such that $\|x\|\ge K$.
\end{ppn}

The previous proposition implies that the distribution of the Markov process $X_t$
converges to its stationary distribution exponentially fast.
More precisely, using  the results e.g. in \cite{Hairer10convergenceof},

\begin{cor}\label{cor:exponential_convergence}
Assume $\ncond$, and let $\alpha >0$ be as in Proposition \ref{ppn:geometric_bound}.
Then there exist constants $C >0$ and $0<\rho < 1$  such that, for all $x$,
$$ \| \BB P_x[X_t=\cdot] - \pi(\cdot) \|_{\tn{TV}} \le C \rho^t e^{\alpha \|x\|}.$$
\end{cor}

In the previous corollary, $\BB P_x$ is the probability under which $\BB P_x[X_0=x]=1$, and
\begin{align}\nn
    \|P-P'\|_{\textnormal{TV}}=\sup_{A\subset \CAL X}|P(A)-P'(A)|
\end{align}
stands for the total variation distance between the probabilities $P$ and $P'$.

\subsection{Bounds}\label{subsection:bounds}

To be able to get  bounds on the variance of the process, we also need to 
find a lower bound for the process.
We can do that by constructing a simple stochastically  smaller process (queue by queue).
The dynamics of the lower bound process for queue $i$ are as follows.
\begin{itemize}
\item
As soon as an item of class j that could be matched to class i items arrives
in the system, one class-$ i $ customers leaves the system.

\item
In case of impatience, each item present in the queue leaves with rate $\gamma(i)$.
\end{itemize}
It is easy to prove by coupling that such a system is a stochastic bound for queue $i$.
Moreover, the $ i $-th marginal of this  process is just a M/M/1 system with arrival rate $\lambda(i)$
and service rate $s_i(x(i))=\sum_{i \sim j} \lambda(j) + \gamma(i) x(i)$.
Hence the marginal of the stationary measure $ \tilde \pi $ is
$$\int \bb1\{ x(i)=z\} \tilde \pi(\dd x)\propto \prod_{z'=0}^{z-1} \frac{\lambda(i)}{ s_i(z'+1)}, \quad z\in\N_0,$$ 
being the empty product defined as $ 1 $.
In the case $\gamma(i)=0, \ \forall i$, we obtain
$$\int x(i) \tilde \pi(\dd x)
=  \frac{\tilde \rho_i}{1- \tilde \rho_i}  , \text{  with  }  \tilde \rho_i =  \frac{\lambda(i)}{\lambda(E(i))}  .$$
Then
\begin{align}\nn
\frac{\tilde \rho_i}{1- \tilde \rho_i}\le 
\int x(i)  \pi(\dd x)
\le 
\int \|x\|_\infty  \pi(\dd x)
\end{align}
for every $ i $, so
\begin{ppn}\label{ppn:lower}
	Under $ \ncond $,
\begin{align}\label{lower_bound}
\int \|x\|_\infty  \pi(\dd x)
\ge
\max_{i\in V}\frac{\tilde \rho_i}{1- \tilde \rho_i},
\end{align}
being $ \pi $ the unique invariant distribution.
\end{ppn}

As a by-product of the previous drift inequalities, we can also deduce interesting bounds on the stationary means of the process. 
For simplicity, we restrict in this subsection to the case
without reneging and measurement errors supported in a bounded interval,
but the same method may be explored to get bounds in the general case.

\begin{cor}\label{cor-1st_moment}
	Assume that $R=\emptyset$ and that there exists  $B\ge 0$ such that $\BB P[-B\le U_{(j,i)}\le B]=1$ for every $(i,j)\in V^2$ such that $\la i,j\ra \in E$.
	Suppose that $ \lambda $ satisfies $ \ncond $, and let $\pi$ be the unique invariant distribution whose existence is guaranteed by Theorem \ref{thm:stability}. 
	Then
	\begin{align}\label{inequality-416}
	\lim_{t\to\infty}\BB E[\|X_t\|_\infty]=
	\int \|x\|_\infty  \pi(\dd x)
	\le 
	\frac{\lambda(V)}{\eta}\Big(\frac{1}{2}+(\check w+2B)|V|\Big).
	\end{align}
\end{cor}

\begin{ppn}\label{ppn-cubic}
	Let $ f_3:\CAL X\to\BB R $ be the (cubic Lyapunov) function defined as $ f_3(x)=\sum_{i\in V}x(i)^3 $.
	Then
	\begin{align}\label{inequality-605}
	\L f_3(x)\le \lambda(V)
	+6\|x\|_\infty[(\check w+2B)|V|+\lambda(V)]-3\eta\|x\|^2_\infty
	\end{align}
	for every $ x\in\CAL X $.
\end{ppn}
	
	Integrating \eqref{inequality-605} with respect to the invariant distribution and using \eqref{inequality-416}, we obtain an upper bound for the second moment $ \int \|x\|_\infty^2 \pi(\dd x) $.
	Put this bound with \eqref{lower_bound} together to obtain the following
	
	\begin{cor}\label{cor:variance}
		Under $ \ncond $,
	\begin{align}\label{inequality-944}
	\lim_{t\to\infty}\BB V(\|X_t\|_\infty^2)
	&=\int \|x\|_\infty^2 \pi(\dd x)
	-\Big[\int \|x\|_\infty \pi(\dd x)\Big]^2
	\\ \nn
	&\le
	\frac{\lambda(V)}{\eta}\Big(
	\frac{1}{3}
	+\frac{1}{\eta}
	[1+2(\check w+2B)|V|]
	[\lambda(V)+(\check w+2B)|V|]
	\Big)
	-\big[\max_{i\in V}\frac{\tilde \rho_i}{1- \tilde \rho_i}\big]^2,
	\end{align}
	being $ \BB V $ the variance associated to $ \BB P $.
	\end{cor}

\section{Simulations}\label{section:simulations}

In this Section we present a few simulations results. In the first set of simulations, we illustrate the 
maximal stability of the Max-Weight policy by representing the long-run behavior of the system. Max-Weight is compared to policies of the {\em priority} type, namely, policies that only take into consideration the rewards of the various matches, but not the queue sizes. 
Then, we give examples of the theoretical bounds obtained in Subsection \ref{subsection:bounds}. Note that these bounds are not meant to be very tight and we leave for future research to optimize
such bounds in particular by better adapting our set of drift inequalities
to the parameters of the model.

\subsection{Comparison with a priority policies.} \label{subsection:greedy}

\begin{figure}[b]
	\begin{subfigure}[b]{0.5\textwidth}
		\centering
		\includegraphics[width=\linewidth]{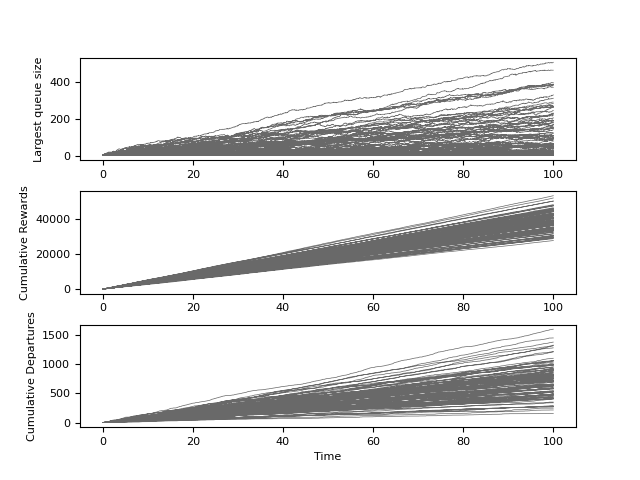}
		\caption*{Priority}
		\label{fig:gull2}
	\end{subfigure}%
	\begin{subfigure}[b]{0.5\textwidth}
		\centering
		\includegraphics[width=\linewidth]{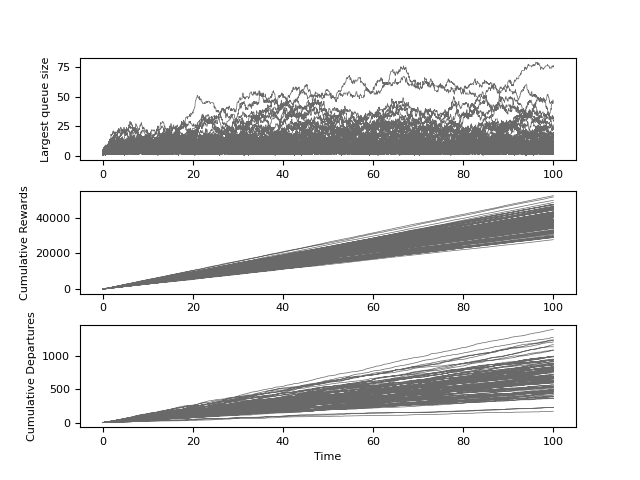}
		\caption*{Max-Weight}
		\label{fig:gull}
	\end{subfigure}%
	\caption{From top to bottom, the graphics represent
		the size of the largest queue, 
		the \\ cumulative reward and the cumulative amount of departures, in the time interval $ [0,100] $.}\label{figure01}
\end{figure}

\begin{figure}[b]
	\centering
	\includegraphics[trim=0pt 0pt 0pt 20pt,width=0.4\textwidth]{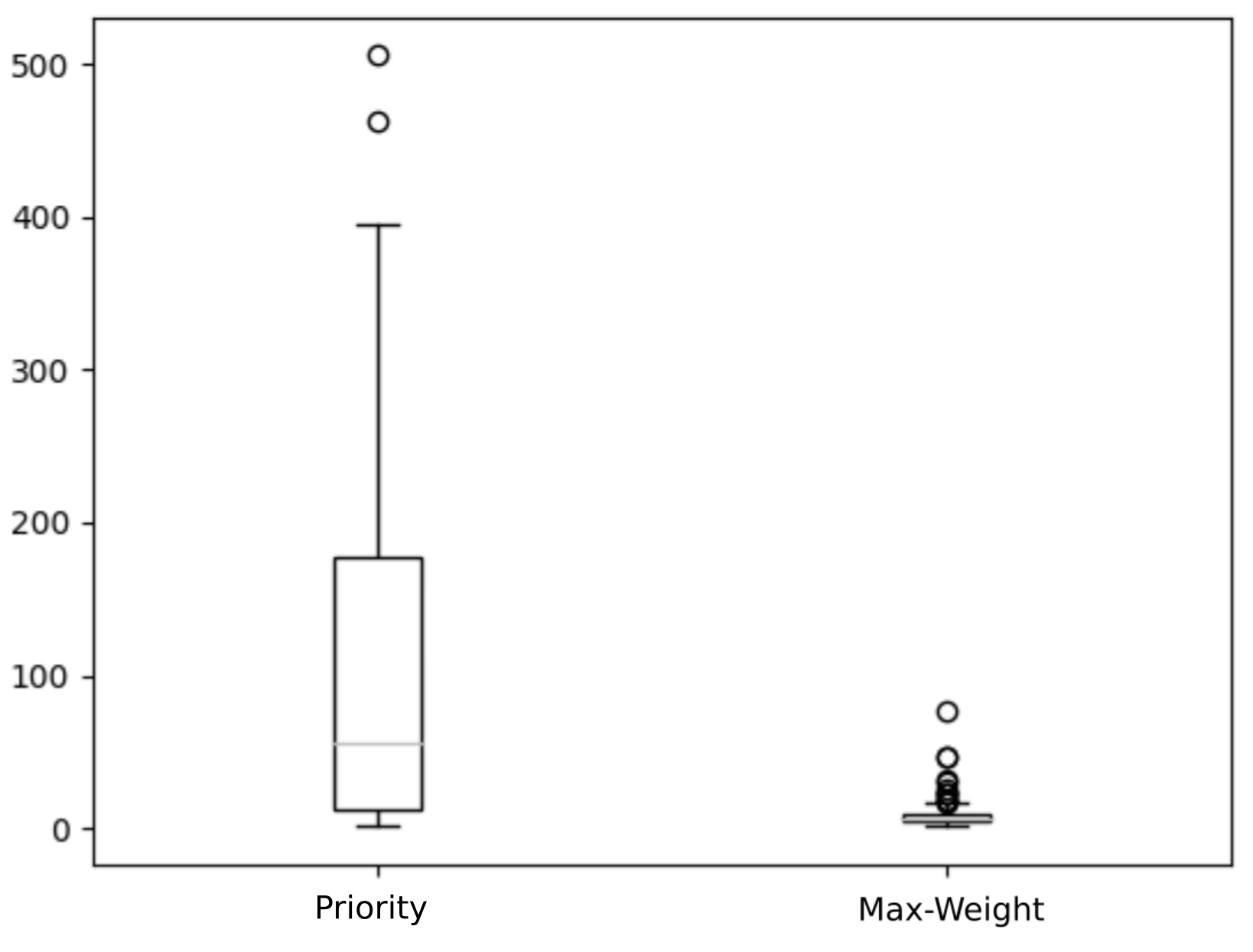}
	\protect\centering
	\caption{ 
		Box-plots for 100 iteration at time $ t=100 $ of the maximum queue size.	} \label{figure02}
\end{figure}

In this subsection, we assume $ F_{(j,i)}=\bb 1_{\R_+} $ for every $ i,j\in V $ such that $ \la i,j\ra \in E $,
namely no measurement errors are considered. A policy is said to be of the {\em priority} type if, whenever the system is in state $ x\in \CAL X $ and a $j$-item enters the system, 
the match of the incoming $j$-item (if any) is chosen uniformly at random from the set
\begin{align}\nn
\argmax_{i\in E(j)\,:\, x(i)>0} w_{( j,i)}.
\end{align}
This amounts to saying that each class-$j$ item has a fixed order of priorities between the neighboring nodes of $i$, 
for choosing its match. 
(Observe that a priority policy can be seen as a Max-Weight policy where we set the errors to $U_{n,(i,j)}\equiv -\infty$ a.s. for all $n,i,j$.) 
As is illustrated in the two examples in Section 5 of \cite{MaiMoy16}, priority policies can be maximal stable or not. On the other hand, as is shown in Theorem 3 of \cite{MoyPer17}, for any graph $G$ outside a particular class of graphs, there always exist a non-maximal priority policy.

To compare the asymptotic behavior of Max-Weight and priority policies for a variety of graphs, arrival and departure rates, and rewards, 
we randomize all these parameters. We perform 100 simulations. % comparing the Max-Weight policy with the y one.
In each one, we first sample the graph $G$, the arrival rates $\lambda $, departure rates $ \gamma $, and rewards $ w $, as follows:
\begin{itemize}
	\item the graph $ G=(V,E) $ is an Erd\H{o}s-R\'enyi one with parameters $ |V|=30 $ and $ p=0.1 $;
	\item the arrival rates $ (\lambda(i))_{i\in V} $ are i.i.d. with common distribution $ \textnormal{Unif}[0,10] $;
	\item the departure rates are i.i.d. with common distribution $ \tn{Ber}(0.5) $ (in average, half of the nodes do not renege);
	\item we consider symmetric rewards, namely $ w_{j,i}=w_{i,j} $ for every $ i,j $ such that $ \la i,j\ra\in E $,  
	and the family $ \{w_{\la i,j\ra}:\la i,j\ra\in E\} $ is i.i.d. with common law
	$ \textnormal{Unif}[0,10] $.
\end{itemize}
In each case, the condition $ \ncond $ is forced, that is, if $\lambda$ is not an element of {\sc Ncond}$(G)$, 
the simulation is re-ran. Then, for each simulation we compare a system ruled by Max-Weight to the same system ruled by the corresponding 
priority policy, for the same sampled parameters. In particular, from Theorem \ref{thm:stability}, the Max-Weight system is necessarily stable, whereas in view of the above observations, the system ruled by the priority policy is not necessarily so.
%The next step is to compare the two policies with these parameters as frozen variables.
%Different simulations correspond to different realizations of the tuple $ (G,\lambda,\gamma,w) $, whose distribution we describe below:

Figure \ref{figure01} represents, in vertical order,
the size of the largest queue,
the cumulative reward (each time a $ (j,i) $-match occurs, a reward $ w_{j,i} $ accumulates)
and the cumulative amount of departures,
the three quantities versus time, in the time interval $ [0,100] $.
Figure \ref{figure02} is the box-plot for the size of the largest queue at time $ t=100 $.

These results put in evidence an important adn striking qualitative difference between the two policies, in the behavior of the buffer size process. Specifically, Max-Weight keeps much less items in the system, whereas both policies have essentially the same performance in terms of cumulative reward and departures. In many cases, the priority policy seems to fail in stabilizing the system, whereas  
the corresponding Max-Weight system is always stable
%even clearly stabilizes the system in many cases, whereas Greedy seems to achieve stability much less often. 
%making in many cases the difference between stability and instability.

\subsection{Bounds in the stationary regime}

\begin{figure}[t]
	\centering
	\includegraphics[width=0.4\textwidth]{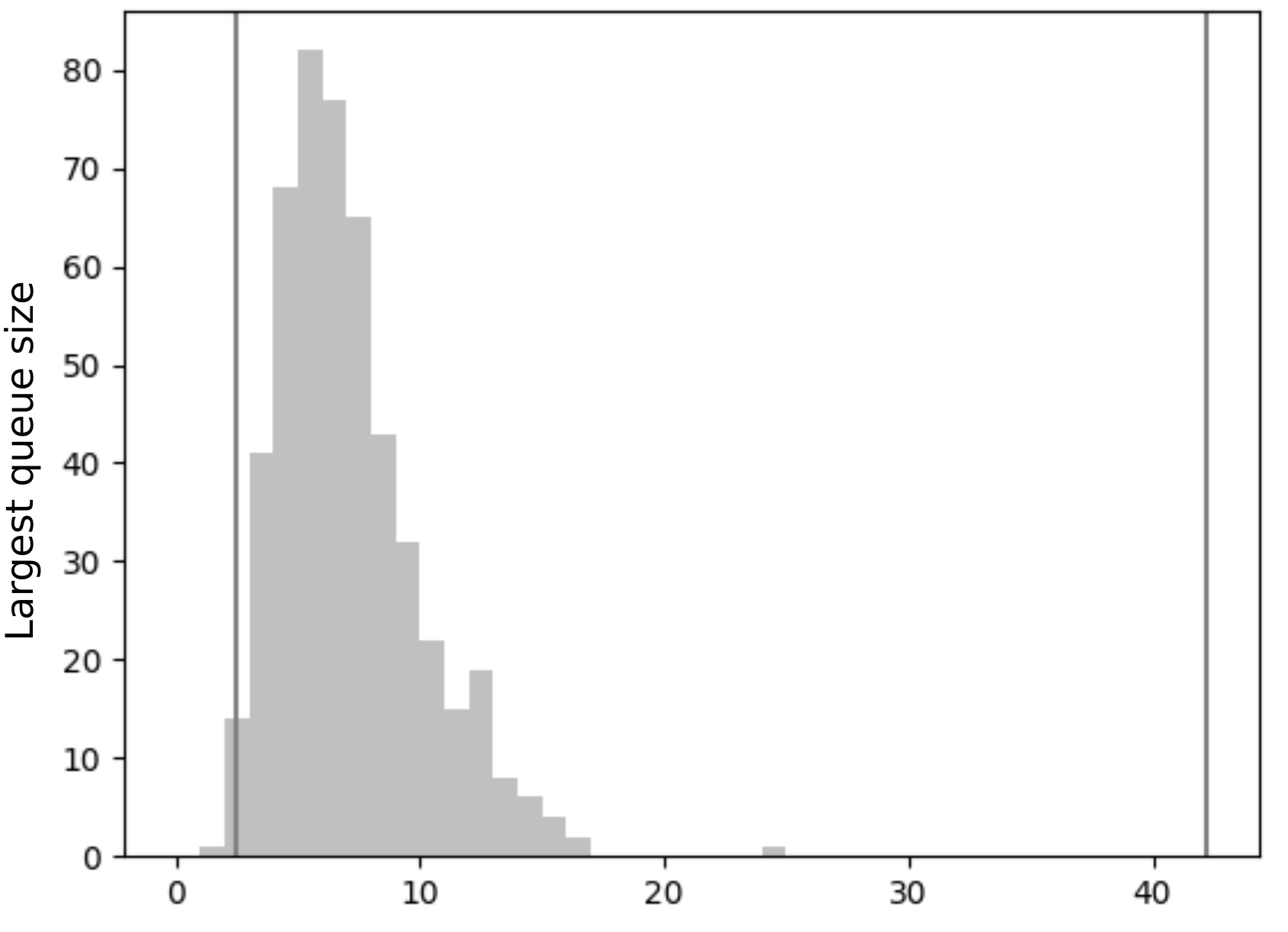}
	\caption{Histogram for 500 simulations at time $ t=100 $ for Match the Longest (i.e. MW with $w_{ij}=0$).} \label{figure03}
\end{figure}

\begin{figure}[t]
	\centering
	\includegraphics[width=0.4\textwidth]{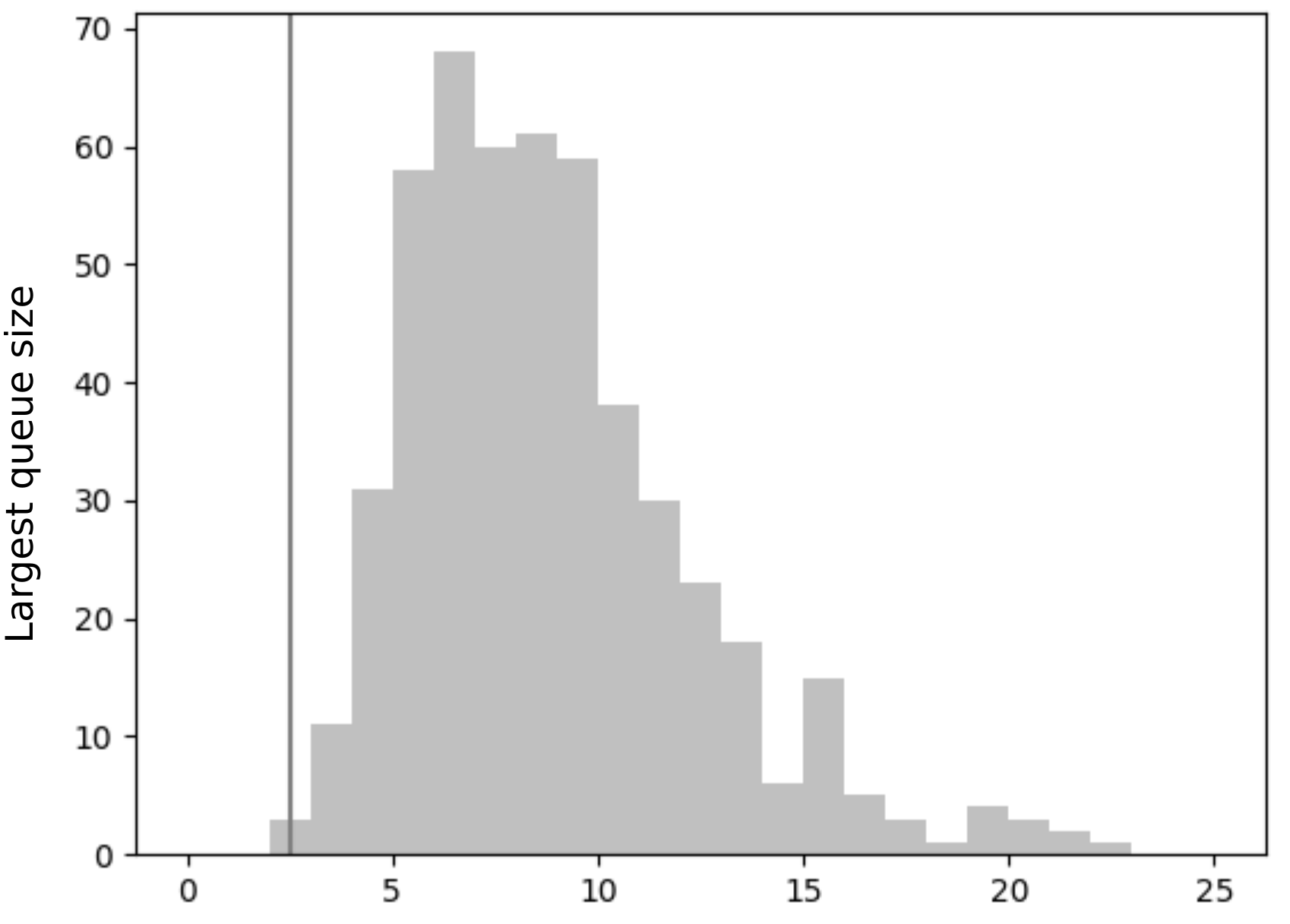}
	\caption{Histogram for 500 simulations at time $ t=100 $ for Max-Weight with rewards and noise} \label{figure04}
\end{figure}

%This subsection is devoted to study how well the theoretical bounds obtained in Subsection \ref{subsection:bounds} adapt to the simulations.
Unlike the previous subsection, impatience is excluded here.

We first restrict our attention to the classical Match the Longest policy, that is, we take $ w_{j,i}=0 $ for every $ i,j\in V $ such that $ \la i,j\ra\in E $,
and we work without measurement errors.
An important difference with respect to the previous subsections is that here the parameters $ (G,\lambda) $ are drawn only once at the beginning.
We clarify that no comparison with a priority policy is made here.
The distribution of the pair $ (G,\lambda) $ is maintained:
Erd\H{o}s-R\'enyi distribution with parameters $ |V|=30 $ and $ p=0.1 $ for $ G $,
and common distribution $ \textnormal{Unif}[0,10] $ for the i.i.d. family $ (\lambda(i))_{i\in V} $.
Figure \ref{figure03} represents a histogram at time $ t=100 $ of 500 simulations.
The vertical lines represent the lower and upper bounds for the expected size of the largest queue respectively given by Proposition \ref{ppn:lower} and Corollary \ref{cor-1st_moment}.
The numerical value for the upper bound for the variance given in Corollary \ref{cor:variance} is $ 7159.85 $.

Finally, we study the case with measurement errors and rewards.
The random variables $ \{U_{j,i}:i,j\in V,\la i,j\ra\in E\} $ giving place to the noises are chosen to be i.i.d. with common distribution $ \tn{Unif}[-1,1] $, and the rewards are again symmetric and i.i.d. with common distribution $ \tn{Unif}[0,10] $.
The distributions of $ G $ and $ \lambda $ do not change,
and the tuple $ (G,\lambda,U,w) $ is sorted only once at the beginning.
Figure \ref{figure04} is a histogram also for 500 iterations at time $ t=100 $, and the vertical line represents the lower bound for the expectation.
The numerical values for the upper bounds for the expectation and the variance are $ 26867.96 $ and $ 13026437.77 $ respectively.

For Match the longest, we observe tight bounds for the expected value,
and weaker bounds for the variance.
All the upper bounds become very bad when we include the presence of rewards or noises.
This is a general picture and not only a particular situation for this choice of the parameters. In the presence of asymetric rewards and noise, the bounds
are not informative. We leave for  future work to obtain better
tailored made bounds in this case.

A typical histogram is very close to the lower bound for the expectation,
and an improvement of the upper bounds is observed if larger values of $ p $ in the distribution of $ G $ are considered.

\section{Proofs}\label{section:proofs}

\subsection{Proof of Proposition \ref{prop:quadratic}}
Fix $x\in\X$ and $\kappa\in (0,\lambda(V))$. It then follows from (\ref{eq:genX}) that 
\begin{align}
\nn
\L f_2(x)
=&\sum_{i  \notin  E(\supp(x))}\lambda(i)\left(2x(i)+1\right)
+\sum_{i  \in  \supp(x)}\Big(\gamma(i)x(i)+\sum_{j  \in   E(i)}\lambda(j)\nu_{x,j}(i) \Big) \left(-2x(i)+1\right)
\\  \nn
=&\lambda(V) + \sum_{i\in \supp(x)}\gamma(i)x(i)\left(-2x(i)+1\right)
+2\sum_{i  \in \supp(x)}x(i)
\Big(\lambda(i) -\sum_{j  \in   E(i)}\lambda(j)\nu_{x,j}(i)\Big)
\\  \nn
=&\lambda(V) + 2\sum_{i\in \supp(x)}x(i)
\Big[\gamma(i)\Big(-x(i)+\frac{1}{2}\Big)+\Big(\lambda(i) -\sum_{j  \in   E(i)}\lambda(j)\nu_{x,j}(i)\Big)\Big]
\\  \nn
=&\lambda(V) + 2\sum_{i\in R}x(i)
\Big[\gamma(i)\left(-x(i)+\frac{1}{2}\right)+\Big(\lambda(i) -\sum_{j  \in   E(i)}\lambda(j)\nu_{x,j}(i)\Big)
\Big]
\\  \nn
&+2\sum_{i\in R^c}x(i)\Big(\lambda(i) -\sum_{j  \in   E(i)}\lambda(j)\nu_{x,j}(i)\Big)
\\ \label{eq:majoregen}
\le& \lambda(V) + 2 \sum_{i\in R}\left[-\gamma(i)x(i)^2+\left(\frac{1}{2} \gamma(i) +\lambda(i)\right)x(i)\right]
+2\sum_{i\in R^c}x(i)\Big(\lambda(i) -\sum_{j  \in   E(i)}\lambda(j)\nu_{x,j}(i)\Big)
\\ \nn
:=&\lambda(V) + A(x)  +B(x),
\end{align}
so the proof is complete in the case $R^c=\emptyset.$ 
Else, we clearly have that 
\begin{equation}\label{eq:majoregen21}
B(x) \le 2\sum_{i\in R^c}x(i)\Big[\lambda(i) - \sum_{j  \in   E(i)}\lambda(j)\nu_{x,j}(i|\maB_\kappa)\pr{\maB_\kappa}\Big]
=:2\sum_{i\in R^c}x(i)r(i),
\end{equation}
where we denote for all $n$ by $\nu_{x,j}(i|\maB_\kappa)$, the probability that a matching $\la i,j \ra$ is performed at $T_n$, given that the state of the system is $x$, the arrival is of class $j$ and 
the event $\maB_{\kappa,n}$ occurs. (Notice that from the iid assumptions, this quantity is independent of $n$.)

Until further notice, we work at an arrival $T_n$, conditionally on the event $\maB_{\kappa,n}$, on the state being $x\in\X$ and the arrival being of class $j$.
We define an auxiliary graph $  G_x=(V_x,E_x) $ by
\begin{align}\nn
V_x=\supp(x)\cap R^c\quad\mbox{ and }\quad
\langle i,j\rangle\in E_x\iff
|x(i)-x(j)|\le 2u_\kappa+\check w,
\end{align}
and we call $  C_1,\ldots, C_L\subset V_x $ its connected components.
We suppose that they are ordered in a decreasing way with respect to the number of items in system, namely, 
if  $ l_1< l_2 $, and if $ i_1\in C_{l_1} $ and $ i_2\in C_{l_2} $,
then $ x(i_1)>x(i_2) $.
The idea behind these definitions is to distinguish nodes whose weights are too far apart.
In particular, we will use the following property: for all $n$, on the event $\maB_{\kappa,n}$, 
if $ C^*=C_1\cup\ldots\cup C_l $ is the union of the first $ l $ clusters and $ C_*=C_{l+1}\cup\ldots\cup C_L $ the union of the remaining ones, 
and we take $ i^*\in C^* $ and $ i_*\in C_* $,
and a $j$-item, for $ j \in E(C^*)\cap E(C_*)$, arrives in the system at $T_n$, then in {\sc mw} this $j$-item prioritizes 
a $i^*$-item over a $ i_* $-item, since 
%\newpage
\begin{align}\nn
[x(i^*)+U_{n,(j,i^*)}]^++w_{( j,i^*)}=x(i^*)+U_{n,(j,i^*)}+w_{( j,i^*)}
&>x(i_*)+U_{n,( j,i^*)}+2u_\kappa+\check w+w_{( j,i^*)}\\
&\ge
[x(i_*)+U_{n,( j,i_*)}]^++w_{( j,i_*)}.\label{property_clusters}
\end{align}
Now, for any $ l\in \llbracket 1,L\rrbracket $ denote $ \check x_l = \max_{i\in C_l}x(i)$.
Observing that 
$$ \check x_l-(2u_\kappa+\check w)|V| \le x(i)\le \check x_l $$
for any $ i\in C_l$,
we get 
\begin{align}
\sum_{i\in C_l}
x(i)r(i)
=\sum_{i\in C_l}^+
x(i)r(i)
+
\sum_{i\in C_l}^-
x(i)r(i)
&\le\check x_l\sum_{i\in C_l}^+
r(i)
+
[\check x_l-(2u_\kappa+\check w)|V|]
\sum_{i\in C_l}^-
r(i)\nn
\\ 
&
=
\check x_l\sum_{i\in C_l}
r(i)
-(2u_\kappa+\check w)|V|
\sum_{i\in C_l}^-
r(i),\label{eq:truc}
\end{align}
where the superscripts $ + $ (resp. $ - $) in the sums mean that we are summing over the indexes $i$'s for which $ r(i)\ge 0 $ (resp. $ r(i)< 0 $). 
But observe that  
	\begin{align*}
	-r(i)=\sum_{j\in E(i)}\lambda(j)\nu_{x,j}(i|\maB_\kappa)\BB P(\maB_\kappa)-\lambda(i)
	\le \sum_{j\in E(i)}\lambda(j)\nu_{x,j}(i)
	\end{align*}
	for all $i \in C_l$,
	implying that 
	\begin{align*}
	-\sum_{l=1}^L\sum_{i\in C_l}^-r(i)
	\le\sum_{l=1}^L\sum_{i\in C_l}^-\sum_{j\in E(i)}\lambda(j)\nu_{x,j}(i)
	&\le\sum_{i\in \supp(x)\cap R^c} \sum_{j\in E(i)}\lambda(j)\nu_{x,j}(i)\\
	&=\sum_{j\in E(\supp(x)\cap R^c)}  \lambda(j)  \sum_{i\in \supp(x)\cap R^c\cap E(j)}\nu_{x,j}(i)\\
	&=\sum_{j\in E(\supp(x)\cap R^c)}  \lambda(j)  
	\le \lambda(V).
	\end{align*}
Therefore, denoting 
\begin{align*}\nn
R_l= 
\sum_{i\in C_l}
r(i),\,l\in \llbracket 1,L \rrbracket, 
\end{align*}
from (\ref{eq:truc}) we readily obtain that 
\begin{equation*}
\sum_{i\in R^c}x(i)r(i) 
= \sum_{l=1}^L \sum_{i\in C_l} x(i)r(i) 
\le \sum_{l=1}^L\check x_l R_l+\lambda(V)(2u_\kappa+\check w)|V|%\label{inequality-58}
\end{equation*}
which, together with \eqref{eq:majoregen21}, yields to 
\begin{equation*}
B(x) \le 2\bigg(\sum_{l=1}^L\check x_l R_l+\lambda(V)(2u_\kappa+\check w)|V|\bigg).
\end{equation*}
%\begin{align}
%\Deltaf_2(x)\le 1+2\sum_{l=1}^L \check x_l R_l+2(\check w+2B)|V|,
%\end{align}
So  we can conclude if we prove that
\begin{align}\nn
\sum_{l=1}^L \check x_l R_l\le (\kappa-\eta)\|x\|_{\infty},
\end{align}
and for this we use the following iterative procedure. 
First, recalling (\ref{eq:tail}), observe that 
\begin{align}
R_1 &= \sum_{i\in C_1}\lambda(i)-\pr{\maB_\kappa}\sum_{j\in E(C_1)}\lambda(j)\sum_{i\in C_1}\nu_{x,j}(i|\maB_\kappa)\nn\\
&= \lambda(C_1)-\pr{\maB_\kappa}\sum_{j\in E(C_1)}\lambda(j)
\le \lambda(C_1)-\left(1-\frac{\kappa}{\lambda(V)} \right)\lambda\left(E(C_1)\right) \le \kappa-\eta,\label{eq:majoreR_1}
\end{align}
where we use, in the last inequality, the fact that $C_1 \in\BB I(R^c)$, %as any non-empty subset of an independent set is also an independent set. 
and in the second identity the following key observation that follows from \eqref{property_clusters}: if the class of an incoming item is an element of $E(C_1) $, then 
it necessarily matches with an item of class in $C_1$, because the queues of nodes from the other clusters $ C_2,\ldots,C_L $ have far fewer stored items. 

We are in the following alternative: if it holds that 
\begin{align}\label{originals}
R_l\le 0 \quad \forall l\in\{1,\ldots, L\}
\end{align} 
we are done since, in view of (\ref{eq:majoreR_1}), 
\begin{equation*}%\label{pera0}
\sum_{l=1}^L\check x_l R_l
\le\check x_1 R_1
%&=\check x_1\left(\sum_{i\in C_1}\lambda(i)-\pr{\maB_\kappa}\sum_{i\in C_1}\sum_{j\in E(i)}\lambda(j)\nu_{x,j}(i|\maB_\kappa)\right)\\
%&=\check x_1\left(\sum_{i\in C_1}\lambda_i-\pr{\maB_\kappa}\sum_{j\in E(C_1)}\lambda(j)\sum_{i\in C_1}\nu_{x,j}(i|\maB_\kappa)\right)\\
%&=\check x_1\left(\lambda(C_1)-\pr{\maB_\kappa}\sum_{j\in E(C_1)}\lambda(j)\right)\\
%&\le \check x_1\left(\lambda(C_1)-\left(1-{\eta\over \lambda(V)}\right)\lambda\left(E(C_1)\right)\right)
%&\le  \|x\|_\infty\left(\lambda(C_1)-\lambda\left(E(C_1)\right)+{\eta\over \lambda(V)}\lambda\left(E(C_1)\right) \right)\\
\le  %\|x\|_\infty\left(-2\eta+\eta\right)=
(\kappa-\eta) \|x\|_\infty.
\end{equation*}
If \eqref{originals} does not hold, we define
\begin{align}\nn
L'= \max\{ l\in \llbracket 1,L\rrbracket \,:\,R_l>0 \},
\end{align}
which, from (\ref{eq:majoreR_1}), is necessarily strictly larger than 1. 
Then we have that 
\begin{align}\label{silla}
\sum_{l=1}^L\check x_l R_l
\le
\sum_{l=1}^{L'}\check x_l R_l
\le 
\sum_{l=1}^{L'-2}\check x_{l}R_l+\check x_{L'-1}(R_{L'-1}+R_{L'}),
\end{align}
%(if $ L'=2 $, the sum $ \sum_{l=1}^{L'-2}\check x_{l}R_l$ is defined to be zero).
where sums over empty sets are set to zero. 
If $ R_{L'-1}+R_{L'}\le 0 $,
we bound the last quantity by $ \sum_{l=1}^{L'-2}\check x_{l}R_l $ and restart the procedure, that is,
we look for
\begin{align*}
\max\{ l\in\llbracket 1,L'-2\rrbracket:R_l>0 \},
\end{align*}
and so on. If $R_{L'-1}+R_{L'}> 0 $,
we bound the r.h.s. of \eqref{silla} by
\begin{align}\nn
\sum_{l=1}^{L'-3}\check x_{l}R_l+\check x_{L'-2}(R_{L'-2}+R_{L'-1}+R_{L'}),
\end{align}
separate between the cases 
\begin{align}\nn
R_{L'-2}+R_{L'-1}+R_{L'}\le 0
 \mbox{ and } 
R_{L'-2}+R_{L'-1}+R_{L'}> 0,
\end{align}
and proceed as above. It only remains to address the case where all partial sums are positive, that is, 
$\sum_{l=l_0}^{L'}R_{l}>0 $ for any $ l_0\in\llbracket 2,L'-2 \rrbracket$.
In this case, the final step of the iteration gives
\begin{align}
\sum_{l=1}^{L'}\check x_l R_l
\le \check x_1 \sum_{l=1}^{L'} R_l.\label{eq:majorefinale}
\end{align}
But as in (\ref{eq:majoreR_1}) we have that 
\begin{align*}
\sum_{l=1}^{L'}R_l &= \sum_{i\in \bigcup_{i=1}^{L'}C_i}\lambda(i)-\pr{\maB_\kappa}\sum_{j\in E\left(\bigcup_{i=1}^{L'}C_i \right)}\lambda(j)\sum_{i\in C_1}\nu_{x,j}(i|\maB_\kappa)\nn\\
&= \lambda\bigg(\bigcup_{i=1}^{L'}C_i\bigg)-\pr{\maB_\kappa}
\lambda\bigg[E\bigg(\bigcup_{i=1}^{L'}C_i\bigg)\bigg]\le \kappa-\eta,
%&\le \lambda(C_1)-\left(1-{\eta\over \lambda(V)}\right)\lambda\left(E(C_1)\right) \le -2\eta+\eta = -\eta,\label{eq:majoreR_1}
\end{align*}
%\begin{align}\nn
%\sum_{l=1}^{L'}R_l
%&=\mu( C_1\cup\ldots \cup C_{L'} )
%-\sum_{i\in C_1\cup\ldots \cup C_{L'}}p(x,x-\delta_i)\le -\eta, 
%\\ \nn
%&
%=\mu( C_1\cup\ldots \cup C_{L'} )
%-\mu(E( C_1\cup\ldots \cup C_{L'}) )<-\eta,
%\end{align}
which follows again from the facts that $\bigcup_{i=1}^{L'}C_i$ is an independent set of $\mathbb I(R^c)$, and that 
a stored item of class in $\bigcup_{i=1}^{L'}C_i$ leaves the system if and only the incoming item's class is in $E\left(\bigcup_{i=1}^{L'}C_i\right)$, as from  \eqref{property_clusters}, 
the {\sc mw} policy would never match the latter with an item of class in $\bigcup_{i=L'+1}^{L}$. Plugging this in \eqref{eq:majorefinale} we get that 
$$\sum_{l=1}^{L'}\check x_l R_l \le  (\kappa-\eta) \|x\|_\infty,$$ which concludes the proof. 

\subsection{Proof of Theorem \ref{thm:stability}}
Necessity of the condition $\ncond$ can be shown similarly to Proposition 2 in \cite{MaiMoy16}: 
For any $I\in \I(R^c)$ and any $t\ge 0$, all items of class in $I$ that have departed the system before $t$ have necessarily been matched with an element of $E(I)$ arrived before $t$, thus we have that 
\begin{align*}
X_t(I) := \sum_{i\in I} X_t(i) & \ge \sum_{i\in I} X_0(i) - \sum_{i\in E(I)} X_0(i) + \sum_{i\in I}N(i)_t - \sum_{i\in E(I)}N(i)_t\\
          & \ge \sum_{i\in I} X_0(i) + M^I_t + t\left(\lambda(I)-\lambda(E(I))\right),
          \end{align*}
          where $M^I$ is a martingale w.r.t. to the natural filtration of $N$. So it is routine to check that the above process tends a.s. to infinity as $t \to \infty$ if 
          $\lambda(I)-\lambda(E(I))>0$, and is at best null recurrent if $\lambda(I)-\lambda(E(I))=0$. 
          
We now prove sufficiency. Fix $\varepsilon>0$.  
For any $x\in\X$, an immediate algebra shows that $A(x)\ge -\lambda(V)-\varepsilon$ implies that 
\begin{equation*}
%\label{eq:majoregen1}
\sum_{i\in R}2\gamma(i)\left[x(i)-\frac{1}{2}\left(\frac{1}{2}  +\frac{\lambda(i)}{\gamma(i)} \right)\right]^2 
\le \sum_{i\in R}\frac{\left( \gamma(i) +2\lambda(i)\right)^2 }{ 2 \gamma(i)}+\varepsilon+\lambda(V),
\end{equation*}
or in other words, that the restriction $x|_R$ of $x$ to its coordinates in $R$, belongs to the set $\mathcal K|_R:= \N^{|R|} \cap \mathcal E$, where $\mathcal E$ is a $|R|$-dimensional ellipsoid. 
Clearly, $\mathcal K|_R$ is a finite set.  

Then we are in the following alternative. If $R^c=\emptyset$ (in which case $\ncond(G)=(\R_+^*)^{|V|}$), the above argument,
together with (\ref{eq:majoregen}),
shows 
that, whenever $x=x|_R$ lies outside the finite set $\mathcal K|_R$, $\L f_2(x) < -\varepsilon$, and 
we conclude by applying (the continuous-time version of) the Foster-Lyapunov Theorem to the Lyapunov function $f_2$.

Now, if $R^c\ne \emptyset$ and $R\ne \emptyset$, the set $\mathcal K$ defined by 
\[\mathcal K=\left\{x \in \subset \N_0^{|V|}\,:\,x|_R \in \mathcal K|_R\mbox{ and }x(i)=0,\,i\in R^c\right\}\]
is also finite, and we can define the quantity 
\[M = \sup_{x\in\X} A(x) =\max_{x\in \mathcal K} A(x)<\infty.\]
But observe that, since $\lambda\in\textsc{ncond}(G)$, the quantity $\eta$ defined by (\ref{eq:defeta}) is an element of $(0,\lambda(V)]$, so it follows by setting 
$\kappa:=\frac{\eta }{ 2}$ in Proposition \ref{prop:quadratic} that for all $x\in\X$, 
\[\L f_2(x) \le \lambda(V)+M+2 \lambda(V)\left(2u^{\eta/2}+\check w\right)|V|-\eta\| x \|_{\infty}.\]
Therefore we get that  $\L f_2(x) < -\varepsilon$ for any $x\in \X$ that does not belong to the finite set 
\[\Big\{x\in\X:\|x\|_\infty\le \frac{\varepsilon +\lambda(V)+M +2\lambda(V)(2u^{\eta/2}+\check w)|V|^2}{\eta} \Big\},\]
and we conclude again by applying Foster-Lyapunov Theorem. 

Last, if $R^c\ne \emptyset$ and $R= \emptyset$, we again have 
that $0<\eta\le \lambda(V)$. So by setting again $\kappa:=\frac{\eta }{ 2}$ in Proposition \ref{prop:quadratic} we obtain similarly that 
%for all $x\in\X$, \[\L f_2(x) \le \lambda(V)+2 \lambda(V)\left(2u^{\eta/2}+\check w\right)|V|-\eta\| x \|_{\infty}.\]
%Therefore, we get that  
$\L f_2(x) < -\varepsilon$ for any $x\in \X$ that does not belong to the finite set 
\[\Big\{x\in\X:\|x\|_\infty\le \frac{\varepsilon +\lambda(V)+2\lambda(V)(2u^{\eta/2}+\check w)|V|^2}{2\eta} \Big\},\]
and we conclude as above. 

%
%
%
%
%If we find $ \varepsilon>0 $ and $ F\subset \X $ finite such that
%\begin{align}
%\Delta f_2(x)\le -\varepsilon
%\end{align}
%for every $ x\notin F $,
%we can conclude by applying Foster's theorem.
%For example, $ \varepsilon=1 $ and
%\begin{align}\nn
%F=\Big\{x\in\CAL X:\|x\|_\infty<\frac{\check w+2B}{\eta}|V| \Big\}
%\end{align}
%works.

\subsection{Proof of Proposition \ref{prop:Ncond}}
To expedite the proof, first observe that, for any $G$ and $\gamma$, 
$\lambda\in(\R_+^*)^{|V|}$ is an element of $\textsc{Ncond}(G,\gamma)$ defined by (\ref{eq:Ncond}) if and only if $\mu_\lambda$ defined by (\ref{eq:defmu}) is an element of 
\begin{equation*}
\textsc{Ncond}^\prime(G,\gamma)=\left\{\mu\in \CAL P^*_{V}\,:\,\mu(I) < \mu(E(I)) \mbox{ for any }I\in\BB I(R^c)\right\},
\end{equation*}
being $ \CAL P^*_{V} $ the set of strictly positive probabilities defined on $ V $;
similarly, $\lambda$ belongs to $\textsc{Ncond}(G)$ defined by (\ref{eq:Ncondno}) if and only if $\mu_\lambda$ belongs to
\begin{equation*}
\textsc{Ncond}^\prime(G)=\left\{\mu\in \CAL P^*_{V}\,:\,\mu(I) < \mu(E(I)) \mbox{ for any }I\in\BB I\right\}.
\end{equation*}

First, if $G$ is bipartite and $R=\emptyset$, then Theorem 1 of \cite{MaiMoy16} shows that the set  $\ncond^\prime(G,\gamma)=\break\ncond^\prime(G)$ 
is empty. From the above remark, so is the set  $\textsc{Ncond}(G,\gamma)=\textsc{Ncond}(G)$. 

Regarding the converse, let us first assume that $G$ is non-bipartite. Then, again from Theorem 1 in \cite{MaiMoy16}, the set 
$\textsc{Ncond}^\prime(G)$ is non-empty. But plainly, for any $\gamma$, we have $\textsc{Ncond}^\prime(G)\subset \textsc{Ncond}^\prime(G,\gamma)$
because any element of $\I(R)$ (if any) is also an element of $\I$. This shows that $\textsc{Ncond}^\prime(G,\gamma)$ is non-empty, and so is 
$\textsc{Ncond}(G,\gamma)$. 

It only remains to consider the case where $G$ is bipartite (in which case $\textsc{Ncond}^\prime(G)$ is empty again by Theorem 1 of \cite{MaiMoy16}), 
and $\gamma$ is such that $R\ne \emptyset$. Let $V_1 \cup V_2$ be a bipartition of $V$ into two independent sets. 
Then, it follows from Theorem 4.2 in \cite{BGM13} that there exists a probability measure $\mu\in \CAL P^*_V$ such that 
\begin{equation}
\label{eq:propmu}
\mu(V_1)=\mu(V_2)=\frac{1}{2}\quad \mbox{ and } \quad \mu(I)<\mu(E(I))\mbox{ for all }I\in\I\setminus\{V_1,V_2\}.
\end{equation}
Suppose that $R\cap V_1\ne \emptyset$ (the case where $R\cap V_2 \ne \emptyset$ is symmetric). 
Our idea is to modify the measure $\mu$ by transferring a small enough amount of mass from a node of $V_2$ to a node of $R\cap V_1$, so as to 
belong to $\I(R^c)$. For this, we first set a node $i\in V_2$ and a positive real number $\varepsilon$ as follows: 
\begin{itemize}
\item[(i)] If $V_1 \subset R$, $i$ is arbitrary in $V_2$ and 
\[\varepsilon=\frac{1}{2}\Bigl(\mu(i)\wedge \min\left\{1-\mu(j)\,:\,j \in R\cap V_1\right\}\Bigl).\]
\item[(ii)] Else, $i$ is an arbitrary element of $E\left(R^c \cap V_1\right)$ and we set 
\[\varepsilon =  \frac{1}{2}\biggl(\min\left\{\mu(E(I)) - \mu(I)\,:\,I \in \I(R^c \cap V_1)\mbox{ s.t. }i \in E(I)\right\}\wedge \min\left\{1-\mu(j)\,:\,j \in R\cap V_1\right\}\biggl),\]
 which is strictly positive in view of the fact that $V_1$ is not an element of $\I(R^c \cap V_1)$. 
\end{itemize}  
We then construct the measure $\hat\mu$ as follows: 
\[\begin{cases}
\hat\mu(i) &= \mu(i) - \varepsilon;\\
\hat\mu(j) &=\mu(j) + \varepsilon \mbox{ for some arbitrary } j\in R\cap V_1;\\
\hat\mu(k) &=\mu(k) \mbox{ for all other }k\in V.
\end{cases}\]
First observe that since $V_1$ intersects with $R$ it is not an independent set of 
$\I(R^c)$. Second, if $V_2 \in \I(R^c)$ (which is true if and only if $R \subset V_1$), we get 
\[\hat\mu(V_2) = \mu(V_2) - \varepsilon = \mu(V_1)-\varepsilon = \hat\mu(V_1)-2\varepsilon =\hat\mu(E(V_2))-2\varepsilon<\hat\mu(E(V_2)).\]
Now take an independent set $I\in\I(R^c)\setminus\{V_2\}$ such that $I\cap V_1 =\emptyset$. 
Then we have that 
\begin{equation}
\label{eq:presque}
\hat\mu(E(I))-\hat\mu(I) \ge \mu(E(I)) - \mu(I) >0,
\end{equation}
where the first inequality follows from the fact that $E(I)\subset V_1$, $I\subset V_2$, and the immediate observations that $\hat\mu(A)\ge \mu(A)$ and 
$\hat\mu(B)\le \mu(B)$ for any $A\subset V_1$ and $B\subset V_2$, and the second, from the facts that $\I(R^c)\subset \I$, $I\not\in\{V_1,V_2\}$ and (\ref{eq:propmu}) . 
Last, take an independent set $I\in\I(R^c)\setminus\{V_2\}$ such that $I\cap V_1 \ne \emptyset$. 
(Notice that such an element exists only if $V_1$ is not included in $R$, so we are in case 
(ii) above.) Then we have
\begin{align*}
\hat\mu(E(I))-\hat\mu(I) &=  \hat\mu(E(I)\cap V_1) + \hat\mu(E(I) \cap V_2) - \hat\mu(I\cap V_2) - \hat\mu(I\cap V_1)\\
				    &=  \hat\mu(E(I\cap V_2)) - \hat\mu(I\cap V_2) + \hat\mu(E(I\cap V_1)) - \hat\mu(I\cap V_1)\\
				    &\ge \hat\mu(E(I\cap V_1)) - \hat\mu(I\cap V_1)\\
				    &\ge \mu(E(I\cap V_1))-\varepsilon - \mu(I\cap V_1)>0,
				    \end{align*}
where in the second equality we use the fact that the graph is bipartite, implying that $E(I)\cap V_1 = E(I\cap V_2)$ and $E(I)\cap V_2 = E(I\cap V_1)$, in the 
first inequality we use (\ref{eq:presque}) in the case where $I\cap V_2 \ne \emptyset$, in the second equality we use the fact that $\hat\mu(A)\ge \mu(A)$ and 
 $\hat\mu(B)\ge \mu(B)-\varepsilon$ for any $A\subset V_1$ and $B\subset V_2$, and in the last inequality, the definition of $\varepsilon$ in case (ii) above. 
 
We have thus proven that $\hat\mu(I)<\hat E(I)$ for any $I\in\I(R^c)$. This implies that $\hat\mu$ is an element of $\textsc{Ncond}^\prime(G,\gamma)$, and thereby, 
an element of $\textsc{Ncond}(G,\gamma)$, that let us  conclude. 
%where, in the first inequality we use the facts that $\hat\mu(B) \le \mu(B)$ for any $B\subset V_2$ and that $\hat\mu(A)=\mu(A)$ for any $A\subset V_1 \cap R^c$; in the second inequality we use (\ref{eq:propmu}) and the fact that $I\cap V_1 \ne V_1$ and $I\cap V_2 \ne V_2$, and in the second equality, the fact that the graph is bipartite. 
%Last, 

\subsection{Proof of Proposition \ref{ppn:geometric_bound}}
Neglecting the reneging part, we get
\begin{align}\nn
    \L e^{\alpha\|x\|} 
    & \le
    \sum_{i \notin E(\supp(x))} \lambda(i)  (e^{\alpha \|x+\delta_i\|}-e^{\alpha \|x\|} ) 
    +\sum_{i \in \supp(x)}
    \tau(i)
    (e^{\alpha \|x-\delta_i\|}-e^{\alpha \|x\|} )
\end{align}
with $\tau(i)=\sum_{j\in E(i)}\lambda(j)\nu_{x,j}(i)$.
Taking common factor $e^{\alpha \|x\|}$
and using that $e^y-1\le y+y^2$ if $|y|$ is small enough, 
last expression can be bounded by
\begin{align}\nn
    &e^{\alpha \|x\|}\bigg[
    \sum_{i \notin E(\supp(x))} \lambda(i)  [\alpha (\|x+\delta_i\|-\|x\|)+\alpha^2] 
    +\sum_{i \in \supp(x)}
    (\alpha( \|x-\delta_i\|- \|x\|)+\alpha^2 )
    \tau(i)
    \bigg]
    \\ \nn
    &\quad =
    e^{\alpha \|x\|}
    \alpha
    \bigg[
    \sum_{i \notin E(\supp(x))} \lambda(i)   (\|x+\delta_i\|-\|x\|) 
    +\sum_{i \in \supp(x)}
    ( \|x-\delta_i\|- \|x\|)
    \tau(i)
    +\alpha\lambda(V)
    \bigg].
\end{align}
From a Taylor developement,
$$\|x \pm \delta_i\|-\|x\| \le \frac{ \pm x(i) }{ ||x||}  + \frac{1}{2\|x\|}.$$
Using this inequality and proceeding as in
the proof of Proposition \ref{prop:quadratic},
we get
\begin{align}\nn
    &\sum_{i \notin E(\supp(x))} \lambda(i)   (\|x+\delta_i\|-\|x\|) 
    +\sum_{i \in \supp(x)}
    ( \|x-\delta_i\|- \|x\|)
    \tau(i)
    \\ \nn
    &\quad\le
    \frac{1}{\|x\|}\bigg(\sum_{i\in \supp(x)}x(i)[\lambda(i)-\tau(i)]+\frac{1}{2}\lambda(V)
    \bigg)
    \\ \nn
    &\quad \le
    \frac{1}{\|x\|}\bigg(
    \lambda(V)(2u_\kappa +\check w)|V|+(\kappa-\eta)\|x\|_\infty
    +\frac{1}{2}\lambda(V)
    \bigg)
    \\ \nn
    &\quad \le
    \frac{\lambda(V)(2u_\kappa +\check w)|V|}{\|x\|}+(\kappa-\eta)\sqrt{|V|}
    +\frac{\lambda(V)}{2\|x\|}
\end{align}
for every $\kappa\in (0,\lambda(V))$.
We have obtained
\begin{align} \nn
    \L e^{\alpha\|x\|} 
    \le
    e^{\alpha \|x\|}
    \bigg(
    \frac{\alpha\lambda(V)(2u_\kappa +\check w)|V|}{\|x\|}+\frac{\alpha\lambda(V)}{2\|x\|}
    +\alpha(\kappa-\eta)\sqrt{|V|}
    +\alpha^2 \lambda(V)
    \bigg).
\end{align}
Since $\eta>0$ due to $\ncond$, we can conclude by choosing $\kappa<\eta$
and $\alpha$ such that $\alpha(\kappa-\eta)\sqrt{|V|}
    +\alpha^2 \lambda(V)<0$.

\subsection{Proof of Corollary \ref{cor-1st_moment}}

In this case, inequality \eqref{inequality-225}
takes the simpler form
\begin{align}\nn
\L f_2(x)\le \lambda(V)
+2 \lambda(V)(2B+\check w)|V|
-2\eta\| x \|_{\infty}.
\end{align}
\eqref{inequality-416} follows after integrating with respect to $\pi$ and using that
\begin{align}\nn
\int \L f(x)\pi(\dd x)=0 \quad \mbox{ for every $ f:\CAL X\to \BB R $ good enough.}
\end{align}

\subsection{Proof of Proposition \ref{ppn-cubic}}

Using that $ f_3(x\pm \delta_i)-f_3(x)=\pm 3x(i)^2+3x(i)\pm 1 $, we get
\begin{align}\nn
\Delta f_3(x)
\le
\lambda(V)
+\sum_{i\in\supp(x)}3x(i)^2r(i)
+\sum_{i\in\supp(x)}3x(i)[\lambda(i)+\tau(i)].
\end{align}
On the one hand,
\begin{align}\nn
\sum_{i\in\supp(x)}3x(i)[\lambda(i)+\tau(i)]
\le 3\|x\|_\infty\lambda(V)+3\|x\|_\infty\sum_{i\in\supp(x)}\tau(i)
\le
6\|x\|_\infty\lambda(V);
\end{align}
on the other,
proceeding as in the proof of Proposition \ref{prop:quadratic},
\begin{align}
\sum_{i\in \supp(x)}3x(i)^2r(i)
\le \sum_{l=1}^L 3\check x_l^2 R_l+6\|x\|_\infty (\check w+2B)|V|
\le 
-3\eta\|x\|_\infty^2+6\|x\|_\infty (\check w+2B)|V|.
\end{align}
Combining these inequalities allows to conclude.

\section*{Acknowledgments}

N. Soprano-Loto is partially supported by
PICT 2015-3154, ANR Grant PRC MATCHES (CE40)
and Aristas.

\bibliographystyle{alphaurl}
\bibliography{biblio}

\end{document}